\documentstyle{article}
\newtheorem{th}{Theorem}
\newtheorem{cor}[th]{Corollary}

\def\sep{\;\vrule\;}
\def\proof#1. {\par
                      \ifdim\lastskip<15pt
                      \removelastskip\penalty-200
                      \vskip15pt plus3pt minus3pt
                      \fi
                       {\def\a{#1}
                       \ifx\a\empty
                       {\noindent\bf Proof.}
                       \else
                       {\noindent\bf Proof of #1.}
                       \fi}\enspace}
\def\restr#1{\,\vrule\,\lower1.75ex\hbox{$#1$}}

\def\endproof{\hfill\hspace{-6pt}\rule[-14pt]{6pt}{6pt}
\vskip22pt plus3pt minus 3pt}

\def\be{\begin{equation}}
\def\ee{\end{equation}}
\def\bea{\begin{eqnarray}}
\def\eea{\end{eqnarray}}
\def\bean{\begin{eqnarray*}}
\def\eean{\end{eqnarray*}}

\def\a{\alpha}
\def\b{\beta}
\def\d{\delta}

\def\e{\varepsilon}

\def\F{\Phi}
\def\g{\gamma}

\def\i{\infty}
\def\k{\kappa}

\def\r{\rho}

\def\t{\theta}

\def\c{{\rm cap}}

\def\ov\overline

\def\c{{\rm cap}}

\font\tenopen = cmbx10
\font\sevenopen = cmbx7
\font\fiveopen = cmbx5
\newfam\openfam
\def\open{\fam\openfam\tenopen}
\textfont\openfam = \tenopen
\scriptfont\openfam = \sevenopen
\scriptscriptfont\openfam = \fiveopen

\def\R{{\open R}}

\def\C{{\open C}}

\def\Z{{\cal Z}}
\def\0{{\bf 0}}
\def\A{{\bf A}}
\def\M{{\cal M}}

\title{Limits of zeros of orthogonal polynomials on the circle\footnote{For  Mathematische Nachrichten issue in memory of F.V. Atkinson}} 
\author{Barry Simon\thanks{Supported in part by  NSF grant DMS-0140592}\and Vilmos  Totik\thanks{Supported by NSF grant DMS-0097484 and  by OTKA  T/034323, TS44782}} 

\begin{document} 

\maketitle 

\begin{abstract} We prove that there is a universal measure on the unit circle such that any 
probability measure on the unit disk is the limit distribution of some subsequence of the 
corresponding orthogonal polynomials. This follows from an extension of a result of Alfaro and Vigil 
(which answered a question of P. Tur\'an): namely, for $n<N$, one can freely prescribe the 
$n$-th polynomial and $N-n$ zeros of the $N$-th one. We shall also describe all possible
limit sets of zeros within the unit disk.
\end{abstract}

\section{Results}

Let $D$ be  the open unit disk. We consider Borel measures $d\mu(t)$, $t\in [-\pi,\pi)$ on the unit 
circle (identified with $\R_{/{\rm mod}\,2\pi}$) of infinite support, and for such a measure 
let 
\be 
\F_n(\mu,z)=z^n+\sum_{k=0}^{n-1}c_{n,k}(\mu)z^k\label{fndef}
\ee
be the $n$-th monic orthogonal polynomial.

It is well known that all zeros of $\F_n$ lie in $D$. The main result of this paper is

\begin{th}\label{AV} For $1\le n<N$, let $\F_n$ be a monic polynomial of degree $n$ with zeros in $D$ 
and let there be given $N-n$ points $a_1,\ldots,a_{N-n}$ in $D$. Then there is a measure $\mu$ on 
the unit circle such that $\F_n(\mu)=\F_n$ and $a_i$, $i=1,\ldots,N-n$, are zeros $($with 
multiplicity$)$ of $\F_N(\mu)$.
\end{th}
In short, one can freely prescribe $\F_n(\mu)$, and $N-n$ zeros of $\F_N(\mu)$.

We do not know if $\F_N(\mu)$ is unique, that is, that the other zeros of $\F_N(\mu)$ are uniquely 
determined by $\F_n$ and $a_1,\ldots,a_{N-n}$.

\begin{cor}\label{sequence} Let $0=n_0<n_1<n_2<\cdots$ be a sequence of natural numbers 
and $\{a_j\}_{j=1}^\i$ a sequence in $D$. Then there is a measure $\mu$ such that for each 
$j=1,2,\ldots$, all $a_i$, $n_{j-1}<i\le n_j$, are zeros $($with multiplicity$)$ of $\F_{n_j}(\mu)$.
 \end{cor}

Indeed, by induction, we get from Theorem \ref{AV} measures $\mu_j$ such that 
\[ 
\F_{n_{j-1}}(\mu_{n_{j-1}})=\F_{n_{j-1}}(\mu_{n_j})
\]
and $a_i$, $n_{j-1}<i\le n_j$, are zeros of $\F_{n_j}(\mu_{n_j})$. By a theorem of Geronimus 
\cite{geronimus}, if for an $m$, the $m$-th orthogonal polynomials for two measures coincide,  
then the same happens for the $l$-th polynomials with all $l\le m$. Hence we have $\F_{n_k}
(\mu_{n_k})=\F_{n_k}(\mu_{n_j})$ for all $k\le j$. Since $\F_{n_k}$ determines the moments 
up to order $n_k$, this shows that for each fixed $m$, the $m$-th moment of the measures 
$\mu_{n_k}$, $k=1,2,\ldots$, are eventually constant, so the weak-$^*$ limit, say $\mu$, of these 
measures exists and obeys $\F_{n_k}(\mu_{n_k})=\F_{n_k}(\mu)$ for all $k$, and the corollary follows.

To formulate the next corollary, recall that the normalized counting measure $\nu_n$ of $\F_{n}$ 
is defined as the measure that puts mass $1/n$ to each zero of $\F_n$ (counting multiplicity). 
We say that a subsequence $\{\F_{n_k}\}$ has zero distribution $\nu$ if the normalized 
counting measures $\nu_{n_k}$ converge to $\nu$ in the weak-$^*$ topology.

\begin{cor}\label{universal} There is a universal measure  $\mu$ on the unit circle such 
that if $\nu$ is any probability measure on the closed unit disk $\bar D$, then some 
subsequence $\{\F_{n_k}(\mu)\}$ has zero distribution $\nu$.
\end{cor}

This is immediate from Corollary \ref{sequence} if we prescribe appropriately $n!-(n-1)!=n! 
((n-1)/n)$ zeros for $\F_{n!}$ in the sequence $\F_0(\mu),\F_{1!}(\mu), \F_{2!}(\mu),\ldots$. 
The rest of the proof is standard (take a countable dense subset $\{\mu_j\}$ in the space 
of probability measures on $\bar D$ so that each $\mu_j$ occurs infinitely often in 
this sequence, and select the $n_k!((n_k-1)/n_k)$ zeros for $\F_{n_k!}$ so that their 
distribution converges weak-$^*$ to $\mu_j$ as we let $n_k$ tend to infinity in such a way 
that $\mu_j=\mu_{n_k}$).

It is known (see, e.g., \cite[Lemma 4]{Widom}) that zeros of orthogonal polynomials cluster 
to the support supp$(\mu)$ of the generating measure $\mu$ if the interior of supp$(\mu)$ 
is empty and $\C\setminus {\rm supp}(\mu)$  is connected. This is due to the fact that 
the $n$-th orthogonal polynomial minimizes the $L^2(\mu)$ norm among all monic polynomials 
of degree $n$. This is the case, for example, if $\mu$ is a measure on the unit circle, but 
its support is not the whole circle. The situation changes if the interior of the support 
is not empty or if $\C\setminus {\rm supp}(\mu)$ is not connected. Consider, for example, 
the closed unit disk or the unit circle and the appropriate (area or arc) Lebesgue measure 
on it, in which case $\F_n(z)=z^n$, hence all the zeros are at the origin. It was P.~Tur\'an 
who asked if for measures on the unit circle it is possible that the zeros of the 
orthogonal polynomials cluster to all points of $D$. First, Szabados \cite{Szabados} had 
shown that for any $\e>0$, there exists a measure for which the set of limit points of the zeros 
had area measure $>\pi-\e$, and then Alfaro and Vigil \cite{alfvigil} noticed that a positive 
answer to Tur\'an's problem follows almost immediately from the Szeg\H{o} recurrence 
\be 
\F_{n+1}(z)=z\F_n(z)-\bar{\a}_n\F^*_n(z), \label{recurr}
\ee
where $\F^*_n(z)=\overline{\F_n(1/\bar z)}$, more precisely from the fact that 
$\mu\leftrightarrow \{\a_n\}_{n=0}^\i$  is a one-to-one correspondence between measures  
on the circle and $D^\i$ (Verblunsky's theorem  \cite{verblunsky}). What Alfaro and Vigil 
showed was that given $\F_n$, one can prescribe one zero of $\F_{n+1}$. Thus, for all $n\ge 1$, 
one can prescribe exactly one zero $a_n$ of $\F_n$, and this is another one-to-one correspondence 
$\mu\leftrightarrow \{a_n\}_{n=1}^\i$ between measures  on the circle and $D^\i$. Indeed, if 
$a_{n+1}$ is a zero of $\F_{n+1}$, then (\ref{recurr}) gives $\overline{\a_n}=a_{n+1} 
\F_n(a_{n+1})/\F^*_n(a_{n+1})$, and all one needs to know is that $|\F_n(z)/\F_n^*(z)|<1$ 
in $D$ (which follows from the maximum principle for holomorphic functions since on the unit 
circle, $|\F_n(z)|=|\F_n^*(z)|$, and $\F_n^*$ has all its zeros outside the unit circle).

Let $L_\mu $ be the set of limit points of the zeros of all $\F_n(\mu)$, $n=1,2,\ldots$. This 
is a subset of $\bar D$, and Tur\'an's problem was if it is possible to have $L_\mu = 
\bar D$. Our second theorem describes the possible sets for $L_\mu $ up to their part 
on the boundary of $D$.

\begin{th}\label{limit} Let $F$ be a compact subset of $\bar D$. Then there is a measure 
$\mu$ on the unit circle such that $L_\mu \cap D=F\cap D$.
\end{th}

What happens on the boundary is less clear. The proof shows that if $F$ contains the unit circle, 
then there is a $\mu$ with $L_\mu =F$. But not every closed subset $F\subseteq \bar D$ is 
an $L_\mu $; for example, $F=[\frac{1}{2},1]$ cannot be the set of limit points. Indeed, suppose 
to the contrary that $[\frac{1}{2},1]=L_\mu $. Then for large $n$, all zeros of $\F_n(\mu)$ lie 
in  the sector ${\rm arg}(z- \frac{1}{4})<\frac{1}{4}$, which gives that $|\F_n(z)|>(\frac{9}{8})^n$ 
for $z$ in a neighborhood $V$ of the point $-1$. Therefore, if $\mu(V)\neq 0$, then 
\[ 
\int|\F_n(\mu)|^2d\mu>c(9/8)^{2n},
\]
which contradicts the fact that the orthogonal polynomials minimize the $L^2(\mu)$-norm among 
all monic polynomials of degree at most $n$ (note that for $z^n$, the $L^2(\mu)$-norm is at 
most the total mass of $\mu$). Thus, we must have $V\cap {\rm supp}(\mu)=\emptyset$, which 
means that the complement of the support of $\mu$ is connected. But then all points in the 
support of $\mu$ are limit points of the zeros (see, e.g., \cite{safftot}), hence the support 
could only consist of the single point 1.

We also note that \cite{Simon} has an alternate way of organizing our proof of Theorem~\ref{AV}.

\smallskip

Derick Atkinson had an important impact on orthogonal polynomials through his book; we are pleased 
to dedicate this paper to his memory.

\section{Proof of Theorem \ref{AV}}

Let $\A {=}(a_1,\ldots,a_{N-n})\in D^{N-n}$, and for $\Z=(z_1,\ldots,z_n)\in D^n$, define 
\be 
d\mu_{\Z,\A}(t)=d_{\Z,\A}\left(\prod_{j=1}^n|e^{it}-z_j|^2\prod_{j=1}^{N-n}
|e^{it}-a_j|^2\right)^{-1} dt, \label{0}
\ee
where $d_{\Z,\A}$ is a normalizing constant for $\mu_{\Z,\A}$ to have total mass $1$. By Geronimus' 
theorem \cite{geronimus}, $\F_N(\mu_{\Z,\A})$ is the monic polynomial with zeros 
\[
z_1,\ldots,z_n,a_1,\ldots,a_{N-n}.
\]
Hence it is enough to show that for some $\Z\in D^n$, the coefficients $c_{n,k}(\mu_{\Z,\A})$, 
$k=0,\ldots,n-1$, from (\ref{fndef}) are the same as $c_k=c_{n,k}(\mu_0)$, where 
\[
\F_n(z)=z^n+\sum_{k=0}^{n-1}c_kz^k
\]
and
\[ 
d\mu_0(t)=\frac{d}{|\F_n(e^{it})|^2}\, dt,
\]
with $d$ a normalizing constant. In other words, we want to show that if 
\[
\M(\mu)=\Bigl(\Re c_0(\mu), \Im c_0(\mu),\Re c_1(\mu), \Im c_1(\mu),\ldots, \Re c_{n-1}(\mu), 
\Im c_{n-1}(\mu)\Bigr),
\] 
then $\M(\mu_{\Z,\A})=\M(\mu_0)$.

Set $m_0=\M(\mu_0)$ and $F_\A(\Z)=\M(\mu_{\Z,\A})$.  We write $z_j=x_j+iy_j$ and we shall 
also consider $F_\A$ as a mapping from $(x_1,y_1,\ldots,x_n,y_n)\in D^n\subset \R^{2n}$ into 
$\R^{2n}$, which gives us a continuously differentiable mapping from an open subset $\R^{2n}$ 
into $\R^{2n}$. With these we want to show that the equation $F_\A(\Z)=m_0$ has a solution in 
$D^n$. We shall do that by showing that the topological degree $d(F_\A,rD^n,m_0)$ of $m_0$ with 
respect to $F_\A$ on $rD^n$ is not zero for $r<1$ sufficiently close to 1, since then 
$F_\A(\Z)=m_0$ has a solution in $rD^n$ (see \cite[Theorem 2.1.1]{lloyd}).

The case $\A=\0$ is instructive, so let us first consider it. Let us also recall (see \cite{lloyd}) 
that $d(F_\0,rD^n,m_0)$ is the sum of the sign of the Jacobian ${\cal J}_{F_\0}$ of $F_\0$ 
(considered as a mapping from $\R^{2n}$ into $\R^{2n}$) at all points $\Z$ which satisfies 
the equation  $F_\0(\Z)=m_0$:
\[ 
d(F_\0,rD^n,m_0)=\sum_{F_\0(\Z)=m_0}{\rm sign}\bigl({\cal J}_{F_\0}\bigr),
\] 
provided none of these $\Z$ is a critical point of $F_\0$.

From Szeg\H{o}'s recurrence (\ref{recurr}), it is immediate that $\F_{n+1}(\mu,z)=z\F_n(\mu,z)$ 
if and only if $z=0$ is a zero of $\F_{n+1}(\mu,z)$, and a repeated application of this gives 
that $z=0$ is an $N-n$ order zero of $\F_N(\mu,z)$ if and only if $\F_N(\mu,z)=z^{N-n}
\F_n(\mu,z)$. In this case, the two associated measures are the same, that is, if 
\[ 
d\mu_{\Z}(t)=d_{\Z}\left(\prod_{j=1}^n|e^{it}-z_j|^2\right)^{-1} dt,
\]
then $F_\0(\Z)=\M(\mu_{\Z,\0})=\M(\mu_\Z)$. Since by Geronimus' theorem, 
\be 
\F_n(\mu_\Z)=\prod_{k}(z-z_k),\qquad \Z=(z_1,\ldots,z_n), \label{fzer} 
\ee
the coefficient $c_{n,k}(\mu_Z)$ equals ($(-1)^{n-k}$ times) the $(n-k)$-th 
elementary symmetric polynomial of the coordinates $z_j$, $j=1,\ldots,n$, of $\Z$. 
These are $n$ analytic functions $g_k(z_1,\ldots,z_n)$, $1\le k\le n$, of the variables 
$z_1,\ldots,z_n$, and in this case, the Jacobian ${\cal J}_{F_0}$ of $F_\0$ is the 
square of the absolute value of the complex Jacobian $|\partial g_k/\partial z_j|_{j,k=1}^n$ 
(see, e.g., \cite[Lemma 2.1]{range}). Thus, the Jacobian  ${\cal J}_{F_\0}$ is everywhere 
nonnegative.

Let $Q\subset \R^{2n}$ be the range of $F_\0$. By Geronimus' theorem, $Q$ is the set consisting 
of (the real and imaginary parts of) the coefficient sequences of all monic polynomials of degree 
$n$ with zeros in $D$, hence $Q$ is a nonempty open subset of $\R^{2n}$. $F_\0$ is the map from 
the zeros of a polynomial to its coefficients, hence it is invariant with respect to permutation 
of the zeros.  Since a polynomial determines its zeros, this map is an $n!$ cover for those 
polynomials with distinct zeros. Furthermore, it is a diffeomorphism in a neighborhood of a 
point of distinct zeros. The set of polynomials with distinct zeros is dense in $Q$, and 
if we select such a $q\in Q$, then it follows that the topological degree $d(F_\0,rD^n,q)$ of 
$q$ with respect to $F_\0$ on $rD^n$ is $n!$ for all $r$ sufficiently close to $1$ (so close 
that $ F_\0^{-1}(\{q\})\subset rD^n$). This is true for a dense set of the $q$'s in $Q$. Therefore, 
we actually have $d(F_\0,rD^n,q)=n!$ for every $q\in Q$ for all $r<1$ sufficiently close to $1$. 
Since $m_0$ is in the range of this mapping ($m_0=F_\0(\Z_0)$ with $\Z_0$ equal to the zero 
sequence of $\F_n$ in some order), we have in particular $d(F_\0,rD^n,m_0)=n!$ for all $r$ 
sufficiently close to $1$ (so close that $\Z_0\in rD^n$).

This has been the case $\A=\0$, and now we turn to general $\A$. Clearly, $F_\A$ is homotopic to 
$F_\0$ under the family $F_{t\A}$, $t\in [0,1]$, and since the degree is invariant under 
homotopy if $m_0$ is not on the images of the boundary, all that is left is to prove there is 
an $r$ sufficiently close to $1$ such that $m_0\notin F_{t\A}(\partial(r D^n))$ for all 
$t\in [0,1]$. In fact, then $d(F_\A,rD^n,m_0)=n!$, and hence $F_\A(\Z)=m_0$ has a solution in 
$rD^n$ (see \cite[Theorem 2.1.1]{lloyd}).

Suppose to the contrary, that for all $k$, there is an $r_k$ with $r_k\to 1$, $\Z_k\in \partial 
({r_k}D^n)$, and $t_k\in [0,1]$ such that $m_0=F_{t_k\A}(\Z_k)$. By selecting a subsequence, we 
may assume that $\Z_k\to \Z^*\in \bar D^n$, $t_k\to t^*$, and $\mu_{\Z_k,t_k\A}\to \mu^*$, 
the latter one in the weak-$^*$ topology. Then at least one component, say $z_1^*$, of $\Z^*$ is 
of absolute value $1$ and each  component $z_{k,j}$ of $\Z_k$ converges to the appropriate 
component $z_j^*$ of $\Z^*$ as $k\to\i$. Hence for the normalizing constants from (\ref{0}), 
we get 
\[
\frac1{d_{\Z_k,t_kA_k}}\ge \int_{-\pi}^\pi\frac{1}{2^{N-1}|e^{it}-z_{k,1}|^2}\, dt=
\frac{2\pi}{2^{N-1}(1-|z_{k,1}|^2)}\to\i
\]
as $k\to\i$. As a consequence, $\mu^*$ is supported on $\{\t\in[-\pi,\pi)\sep e^{i\t}\in 
\Z^*\cap \partial D\}$, that is, its support consists of at most $n$ points, and we may assume 
that $z_1^*\in {\rm supp}(\mu^*)$. If $H_{n-1}$ is any polynomial of degree at most $n-1$, we 
have 
\[
\int \F_n \overline{H_{n-1}}\, d\mu_{\Z_k,t_k\A}=0
\]
for all $k$ (recall that $\F_n$ is the $n$-th orthogonal polynomial with respect to each 
$\mu_{\Z_k,t_k\A}$ by the choice of $\Z_k$). Hence it follows, by weak-$^*$ convergence, that 
\[
\int \F_n \overline{H_{n-1}}\,d\mu^*=0
\] 
for any $H_{n-1}$. Now choose $H_{n-1}$ so that it vanishes at all points of the support of 
$\mu^*$ except for $z_1^*$. Then the previous equality gives 
\[
0\neq \F_n(z_1^*)\overline{H_{n-1}(z_1^*)}\, \mu^*(\{z_1^*\})=
\int \F_n \overline{H_{n-1}} \, d\mu^*=0,
\]
and this contradiction shows that, in fact, $m_0\neq F_{t\A}(\partial rD^n)$ for all 
$t\in [0,1]$, provided $r<1$ is sufficiently close to $1$.\endproof

\section{Proof of Theorem \ref{limit}}

Let $D_r(z)$ denote the disk of radius $r$ about the point $z$.

We shall use a reasoning similar to the one in \cite[Example 2.1.2]{stahltot}. The proof is 
based on the observation that if $\mu$ is a measure consisting of $m$ mass points and a 
very small part somewhere else, then $\F_m(\mu)$ will have precisely one zero close to 
each mass point (and, of course, no other zero), and $\F_{m+1}(\mu)$ will have precisely 
one zero close to each mass point, plus an additional zero, and this additional zero will 
be what will move around in the construction to describe the assumed limit set.

Let $F\subseteq \bar D$ be the given closed set, and select a countable set 
$\{S_n\}_{n=0}^\i$ in $D$ such that $F$ is precisely the set of limit points of 
$\{S_n\}_{n=0}^\i$. Choose also pairwise different points $P_1,P_2,P_3,\ldots$ on the unit 
circle such that $S_n$ lies on the segment joining $P_{2n+1}$ and $P_{2n+2}$. It is 
sufficient to show a measure $\mu$ such that for each $n$, the zeros of $\F_{2n}$ 
lie very close to $\{P_1,P_2,\ldots,P_{2n}\}$, while the zeros of $\F_{2n+1}$ lie very 
close to $\{P_1,P_2,\ldots,P_{2n}\}\cup\{S_n\}$. In fact, then the set of limit points 
of the zeros is the closure of the set $\{P_j\}_{j=1}^\i\cup\{S_n\}_{n=1}^\i$, the 
intersection of which with $D$ is precisely $F\cap D$. The measure $\mu$ will be of the form
\be 
\mu=\sum_{n=0}^\i \e_n\Bigl(\b_n\d_{P_{2n+1}}+(1-\b_n)\d_{P_{2n+2}}\Bigr), \label{mudef}
\ee
where $\d_P$ denotes the Dirac mass at the point $P$ and $\b_n\in(0,1)$, $\e_n>0$ will 
be chosen below.

First of all, we require $\e_{n+1}<\e_n/2$ for all $n$. Start with $\e_0=1$, 
$\b_0= \frac12$. Suppose that we have already selected $\e_0,\ldots,\e_{n-1}$ and 
$\b_0,\ldots,\b_{n-1}$ for some $n>0$. Set 
\[
\pi_n =\prod_{j=1}^{2n}|P_{2n+1}-P_j|, 
\]
\[
\k_n =\prod_{j=1}^{2n}|P_{2n+2}-P_j|, 
\]
and for a $\b\in[0,1]$, minimize the expression 
\[ 
\b\pi_n^2|P_{2n+1}-s|^2+(1-\b)\k_n^2|P_{2n+2}-s|^2
\]
for $s\in\C$. It is clear that the minimum is taken at some point of the segment connecting 
$P_{2n+1}$ and $P_{2n+2}$, and if $\b=0$, then it is taken at $P_{2n+2}$, while if $\b=1$, 
then it is taken at $P_{2n+1}$. As $\b$ moves from $0$ to $1$, there will be a value, which we call 
$\b_n$, for which the minimum is taken at the point $s=S_n$. Now by continuity, there is a $\g_n>0$ 
such that if $|z_j-P_j|<\g_n$ for $j=1,2,\ldots,2n$ and $T_{2n}(z)=\prod_{j=1}^{2n}(z-z_j)$, then 
the minimum of 
\[
\b_{n+1}|T_n(P_{2n+1})|^2|P_{2n+1}-s|^2+(1-\b_{n+1})|T_{2n}(P_{2n+2})|^2|P_{2n+2}-s|^2
\]
is taken somewhere in $D_{1/n}(S_n)$. We may assume that $\g_n$ is smaller than $1/n$ and that it 
is also smaller than half of the minimal distance in the set 
\[
\{P_1,P_2,\ldots,P_{2n+1},P_{2n+2}\}.
\]
Thus,
\begin{description}
  \item[(A)] {\it If $\mu$ is a measure such that $\F_{2n+1}(\mu)$ has a zero
in each $D_{\g_n}(P_j)$, $j=1,\ldots,2n$, then the remaining zero of $\F_{2n+1}(\mu)$ 
will be in $D_{1/n}(S_n)$.} 
\end{description} 
Next, we claim that 
\begin{description}
  \item[(B)] {\it There is an $\eta_n>0$ such that if $\mu$ is of the form 
  \be 
  \mu=\sum_{j=0}^{n-1} \e_j\Bigl(\b_j\d_{P_{2j+1}}+(1-\b_j)\d_{P_{2j+2}}\Bigr)+\nu,
\label{muform}
\ee
where $\nu$ is any  measure supported in $\bar D$ with total mass $\|\nu\|<\eta_n$, 
then $\F_{2n+1}(\mu)$ has a zero in each $D_{\g_n}(P_j)$, $j=1,\ldots,2n$.} 
\end{description} 
In fact, we have 
\[
\int\Bigl|z\prod_{j=1}^{2n}(z-P_j)\Bigr|^2 \,d\mu \le 2^{2n}\|\nu\|\le 2^{2n}\eta_n.
\] 
Furthermore, there is a $\r_n>0$ such that for all polynomials $R_{2n+1}(z)=z^{2n+1}+\cdots$, 
the zeros of which omit at least one of $D_{\g_n}(P_j)$, $j=1,\ldots,2n$, we have 
\[
\max_{j=1,\ldots,2n}|R_{2n+1}(P_j)|>\r_n.
\]
Thus, if
\[
\eta_n<2^{-2n}\r_n^2\, \min_{0\le j\le n-1}\, \e_j\min\{\b_j,1-\b_j\},
\] 
then no such polynomial can minimize the $L^2(\mu)$-norm (recall that $\mu$ has a mass $\ge 
\e_j\min\{\b_j,1-\b_j\}$ at each $P_{2j},P_{2j+1}$, $j=1,\ldots,n$), and since $\F_{2n+1}(\mu)$ 
minimizes the $L^2(\mu)$-norm, the claim follows.

A perfectly similar argument gives that
\begin{description}
  \item[(C)] {\it There is an $\eta_n'>0$ such that if $\mu$ is of the form $(\ref{muform})$ 
  where $\nu$ is any  measure supported in $\bar D$ with total mass $\|\nu\|<\eta_n'$, 
  then $\F_{2n}(\mu)$ has a zero in each $D_{\g_n}(P_j)$, $j=1,\ldots,2n$ $($and, of course, 
  no other zero$)$.} 
\end{description}

Now set $\e_n<\min(\eta_n,\eta_n',\e_{n-1})/2$. With this choice, the measure from (\ref{mudef}) 
satisfies {\bf (B)} and {\bf (C)}, and hence we obtain from {\bf (A)} that we can number the 
zeros $z_1^{(2n+1)},\ldots,z_{2n+1}^{(2n+1)}$ of $\F_{2n+1}(\mu)$ in such a way that for 
$j=1,2,\ldots,2n$, we have $|z_j^{(2n+1)}-P_j|<1/n$, and $|z_{2n+1}^{(2n+1)}-S_{n+1}|<1/n$, 
and similarly, it follows from {\bf (C)} that we can number the zeros $z_1^{(2n)},\ldots, 
z_{2n}^{(2n)}$ of $\F_{2n}(\mu)$ in such a way that for $j=1,\ldots,2n$, we have $|z_j^{(2n)}-P_j| 
<1/n$, and this is what we wanted to achieve. \endproof

\vskip1cm

\noindent Barry Simon

Department of Mathematics

California Institute of Technology

Pasadena, CA 91125

USA

{\it bsimon@caltech.edu}
\bigskip

\noindent Vilmos Totik

Bolyai Institute

University of Szeged

Szeged

Aradi v. tere 1, 6720, Hungary

\smallskip
and
\smallskip

Department of Mathematics

University of South Florida

4202 E. Fowler Ave, PHY 114

Tampa, FL 33620-5700 

USA

{\it totik@math.usf.edu}
\end{document}